\documentclass[11pt]{article}
\usepackage{latexsym}
\usepackage[all]{xy}
\usepackage{amsmath}
\usepackage{amssymb}
\usepackage{amsfonts}

\newtheorem{thm}{Theorem}[section]
\newtheorem{prop}[thm]{Proposition}
\newtheorem{lem}[thm]{Lemma}
\newtheorem{df}[thm]{Definition}
\newtheorem{cor}[thm]{Corollary}

\newtheorem{rmk}[thm]{Remark}
\newtheorem{ex}[thm]{Example}

\begin{document}

\title{\textbf{A remark on $K$-theory and $S$-categories}}
\bigskip
\bigskip

\author{\bigskip\\
Bertrand To\"en \\
\small{Laboratoire Emile Picard}\\
\small{UMR CNRS 5580} \\
\small{Universit\'{e} Paul Sabatier, Toulouse}\\
\small{France}\\
\bigskip \and
\bigskip \\
Gabriele Vezzosi \\
\small{Dipartimento di Matematica}\\
\small{Universit\`a di Bologna}\\
\small{Italy}\\
\bigskip}

\date{}

\maketitle

\begin{abstract}
It is now well known that the $K$-theory of a Waldhausen category depends on more than just
its (triangulated) homotopy category (see \cite{sch}). The purpose of this note is to show that
the $K$-theory spectrum of a (good) Waldhausen category is completely determined
by its Dwyer-Kan simplicial localization, without any additional structure. As the
simplicial localization is a refined version of the homotopy category which also determines the triangulated structure, our result is a possible answer to the general question:
\textit{``To which extent $K$-theory is not an invariant of triangulated derived categories ?''}
\end{abstract}

\medskip

\textsf{Key words:} $K$-theory, simplicial categories, derived categories.

\textsf{MSC-class:} $55$N$15$, $19$D$10$, $18$E$30$, $18$G$55$, $18$D$20$.

\medskip

\tableofcontents

\bigskip

\setcounter{section}{0}

\begin{section}{Introduction}

As recently shown by Marco Schlichting in \cite{sch}, the $K$-theory spectrum (actually the $K$-theory groups)
of a stable model category depends on strictly more than just its triangulated homotopy category; indeed he exhibits two
Waldhausen categories having equivalent (triangulated) homotopy categories and non weakly equivalent associated $K$-theory spectra.
Because of this there is no longer any hope of defining a reasonable $K$-theory functor on the level
of triangulated categories (see \cite[Prop. 2.2]{sch}). In this paper we show that, if one replaces in the above statement
the homotopy category (i.e. the Gabriel-Zisman localization with respect to weak equivalences), with the more refined
\textit{simplicial  localization} of Dwyer and Kan, then one actually gets an invariance statement; more precisely, we
prove that the $K$-theory spectrum of
a \textit{good} Waldhausen category (see Definition \ref{d1}) is an
invariant of its simplicial localization without any
additional structure. As the
simplicial localization is a refined version of the homotopy category,
that is a simplicially enriched category lying in between the category itself and its homotopy category, we like to consider this result
as a possible answer to the general question:
\textit{``To which extend $K$-theory is not an invariant of triangulated derived categories ?''}.
In a sense, our result explains exactly what ``structure''
is lacking in the derived (or homotopy) category of a good Waldhausen category, in order to reconstruct its $K$-theory. \\

Our approach consists first in defining a $K$-theory functor on the level of $S$-categories (i.e. of simplicially enriched
categories) satisfying some natural properties, and then in proving that, when applied to
the simplicial localization of a \textit{good} Waldhausen category $C$, this construction yields
a spectrum which is weakly equivalent to the Waldhausen's $K$-theory spectrum of $C$. \\

\textbf{Good Waldhausen categories.} Let us briefly describe the class of
Waldhausen categories for which our result holds (see Definition \ref{d1} for details and the last paragraph of the Introduction for our conventions and notations on Waldhausen categories).
Roughly speaking, a \textit{good} Waldhausen category is a Waldhausen category that can be embedded
in the category of fibrant objects of a pointed model category, and whose Waldhausen structure
is induced by the ambient model structure (Definition \ref{d1}). Good Waldhausen categories
behave particularly well with respect to simplicial localization as they possess a nice homotopy
calculus of fractions (in the sense of \cite{dk3}). The main property of good Waldhausen categories is the following
form of the approximation theorem.

\begin{prop}\emph{(See Proposition \ref{p1})}\label{p0}
Let $f : C \longrightarrow D$ be an exact functor between good Waldhausen categories.
If the induced morphism $L^{H}C \longrightarrow L^{H}D$ between the simplicial localizations is an equivalence
of $S$-categories, then the induced morphism
$$K(f) : K(C) \longrightarrow K(D)$$
is a weak equivalence of spectra.
\end{prop}

Though there surely exist
non-good Waldhausen categories (see Example \ref{nongood}),
in practice it turns out that given a Waldhausen category there is always a good Waldhausen
model i.e. a good Waldhausen category with the same $K$-theory space up to homotopy. For example, the category of perfect complexes on
a scheme and the category of spaces having a given space as a retract, both have
good Waldhausen models (see Example \ref{ex}); this shows that the class of good Waldhausen categories contains interesting
examples\footnote{Actually, we do not know any reasonable example
for which there is no good Waldhausen model.}. \\

\textbf{$K$-Theory of $S$-categories.} For an $S$-category
$T$, which is \textit{pointed} and has \textit{fibered products} (see definitions \ref{d3} and \ref{d4} for details),
we define an associated good Waldhausen category $M(T)$, by embedding $T$ in the model category
of simplicial presheaves on $T$. The $K$-theory spectrum $K(T)$ is then
defined to be the Waldhausen $K$-theory spectrum of $M(T)$
$$K(T):=K(M(T)).$$
Proposition \ref{p0} immediately implies that $K(T)$ is invariant,
up to weak equivalences of spectra, under equivalences in the argument $T$ (see the end of Section $3$).

If $C$ is a good Waldhausen category, then its simplicial localization $L^{H}C$
is pointed and has fibered products (Proposition \ref{p2}). One can therefore consider its $K$-theory
spectrum $K(L^{H}C)$. The main theorem of this paper is the following

\begin{thm}\label{unodue}\emph{(See Theorem  \ref{t1})}
If $C$ is a good Waldhausen category, there exists a weak equivalence
$$K(C)\simeq K(L^{H}C).$$
\end{thm}

\medskip

As a main corollary, we get the following result that actually motivated this paper

\begin{cor}\label{unotre}\emph{(See Corollary \ref{c1})}
Let $C$ and $D$ be two good Waldhausen categories. If the two $S$-categories
$L^{H}C$ and $L^{H}D$ are equivalent, then the  $K$-theory spectra
$K(C)$ and $K(D)$ are isomorphic in the homotopy category of spectra.
\end{cor}

Another interesting consequence is the following.

\begin{cor}\emph{(See Corollary \ref{c2})}
Let $M_{1}:=m\mathcal{M}(\mathbb{Z}/p^{2})$ and $M_{2}:=m\mathcal{M}(\mathbb{Z}/p[\epsilon])$ be the two stable model categories considered in
\cite{sch}. Then, the two $S$-categories $L^{H}M_{1}$ and $L^{H}M_{2}$ are not
equivalent.
\end{cor}

This last corollary implies the existence of two \textit{stable $S$-categories} (see Section $7$),
namely $L^{H}M_{1}$ and $L^{H}M_{2}$, with equivalent triangulated homotopy categories, but which are not equivalent. \\

\textbf{What we have not done.}
To close this introduction let us mention that we did not investigate the
\textit{full} functoriality of the construction $T \mapsto K(T)$ from $S$-categories to spectra, and more generally we did
not try to fully develop the $K$-theory of $S$-categories, though we think that this deserves to be done in the future. In a similar vein, we think that the equivalence
of our main theorem (Thm. \ref{unodue}) is in a way functorial in $C$, at least up to homotopy, but we did not try to prove this. Thus, the results of this paper definitely do not pretend to be
optimal, as our first motivation was only to give a proof of Corollary \ref{c1}.
However, the interested reader might consult the last section in which we present
some ideas towards more intrisinc constructions and results, independent of the notion of Waldhausen category. \\

\bigskip

\textbf{Organization of the paper.} In Section 2 we introduce the class of Waldhausen categories (\textit{good Waldhausen categories}) for which our main result holds; for such categories $C$, we prove that the geometric realization of the subcategory $W$ of weak equivalences is equivalent to the geometric realization of the $S$-category of homotopy equivalences in the hammock localization $L^{H}(C)$ of $C$ along $W$. We also list some  examples of good Waldhausen categories. In Section 3 we define DK-equivalences and prove Proposition \ref{p0}, a strong form of the approximation theorem for good Waldhausen categories. In Section 4 we define when an $S$-category is pointed and has fibered products, and prove that the hammock localization of a good Waldhausen category along weak equivalences is an $S$-category of this kind. In Section 5 we define the $K$-theory of pointed $S$-categories with fibered products and study its functoriality with respect to equivalences. Section 6 contains the main theorem showing that the $K$-theory of a good Waldhausen category is equivalent to the $K$-theory of its hammock localization along weak equivalences. As a corollary we get that the hammock localization completely determines the $K$-theory of a good Waldhausen category. Finally, in Section 7 we discuss possible future directions and relations with other works.  \\

\medskip

\textbf{Conventions and review of Waldhausen categories, $S$-categories and simplicial localization.}
Throughout this paper, a \textit{Waldausen category} will be the dual of a usual Waldhausen category, i.e. our Waldhausen categories we will always be \textit{categories with fibrations and weak equivalences}
satisfying the axioms \textit{duals} to Cof1, Cof2, Cof3 (\cite[1.1]{wa}),
Weq1 and Weq2 (\cite[1.2]{wa}). The reason for such a choice is only stylistic, in order to avoid having to dualize too many times in the text. \\
Explicitly, in this paper a Waldhausen category will be a triple $(C, \mathrm{fib}(C), \mathrm{w}(C))$ consisting of category $C$ with subcategories $\mathrm{w}(C)\hookrightarrow C$ and $\mathrm{fib}(C)\hookrightarrow C$ whose morphisms will be called (weak) equivalences and fibrations, respectively, satisfying the following axioms
\begin{itemize}
    \item $C$ has a zero object $*$.
    \item $(\mathbf{Cof}1)^{\mathrm{op}}$: The subcategory $\mathrm{fib}(C)$ contains all isomorphisms in $C$.
    \item $(\mathbf{Cof}2)^{\mathrm{op}}$: For any $x\in C$, the morphism $x\rightarrow *$ is in $\mathrm{fib}(C)$.
    \item $(\mathbf{Cof}3)^{\mathrm{op}}$: If $x\rightarrow y$ is a fibration, then, for any morphism $y'\rightarrow y$ in $C$, the pullback $x\times_{y}y'$ exists in $C$ and the canonical morphism $x\times_{y}y'\rightarrow y'$ is again a fibration.
    \item $(\mathbf{Weq}1)^{\mathrm{op}}$: The subcategory $\mathrm{w}(C)$ contains all isomorphisms in $C$.
    \item $(\mathbf{Weq}2)^{\mathrm{op}}$: If in the commutative diagram $$\xymatrix{y \ar[d] \ar[r]^{p} & x \ar[d] & z \ar[l] \ar[d]\\ y' \ar[r]_{p'} & x' & z' \ar[l]}$$ in $C$, $p$ and $p'$ are fibrations and the the vertical arrows are equivalences, then the induced morphism $y\times_{x}z\rightarrow y'\times_{x'}z'$ is again an equivalence.
\end{itemize}

To any usual Waldhausen category there is an associated $K$-theory spectrum (or space) as defined in \cite[\S 1.3]{wa} using Waldhausen $S_{\bullet}$-construction. If $C$ is a Waldhausen category according to our definition above, then there is a dual $S_{\bullet}$-construction, denoted by $S_{\bullet}^{\mathrm{op}}$, formally obtained by replacing cofibrations with fibrations (with opposed arrows) and pushouts with pullbacks in the usual $S_{\bullet}$-construction. The dual $S_{\bullet}$-construction applied to our $C$, produce a $K$-theory spectrum $$n\mapsto \left|wS_{\bullet}^{\mathrm{op}}\cdots S_{\bullet}^{\mathrm{op}}C\right|$$ denoted by $K(C)$. Note that $K(C)$  obviously coincides with the usual Waldhausen $K$-theory spectrum (as defined in \cite[\S1.3 p. 330]{wa}) of the dual category $C^{op}$, considered as a usual Waldhausen category (i.e. a category with equivalences and cofibrations satisfying the dual of the above axioms).\\

For \textit{model categories} we refer to \cite{ho} and \cite{hi} which are standard references on the subject. We will often use a basic link between model categories and Waldhausen categories, namely the fact that if $M$ is a pointed model category (i.e. a model category in which the initial object $\emptyset$ and the final object $*$ are isomorphic), then its subcategory $M^{f}$ of fibrant objects together with the induced subcategories of equivalences and fibrations is in fact a Waldhausen category (according to our convention). This follows from \cite[Thm. 19.4.2 (2) \& Prop. 19.4.4 (2)]{hi}, and can also be checked more elementarily.\\

By an $S$\textit{-category}, we will mean a category enriched over the category of simplicial sets. If $T$ is an $S$-category, we will denote by $\pi_{0}T$ the category with the same objects as $T$ and with morphisms given by $\mathrm{Hom}_{\pi_{0}T}(x,y):=\pi_{0}(\underline{\mathrm{Hom}}_{T}(x,y))$, where  $\underline{\mathrm{Hom}}_{T}(x,y)$ is the simplicial set of morphisms between $x$ and $y$ in $T$.
Recall the following fundamental definition.

\begin{df}\label{eqS}
Let $f : T \longrightarrow T'$ be a morphism of $S$-categories.
\begin{enumerate}
\item The morphism $f$ is \emph{essentially surjective} if the induced functor $\pi_{0}f : \pi_{0}T \longrightarrow
\pi_{0}T'$ is an essentially surjective functor of categories.
\item The morphism is \emph{fully faithful} if for any pair of objects $(x,y)$ in $T$,
the induced morphism $f_{x,y} : \underline{Hom}_{T}(x,y) \longrightarrow \underline{Hom}_{T'}(f(x),f(y))$
is an equivalence of simplicial sets.

\item The morphism $f$ is an \emph{equivalence} if it is essentially surjective and fully faithful.

\end{enumerate}
\end{df}

Given a category $C$ and a subcategory $W$, W. Dwyer and D. Kan have defined
in \cite{dk3} an $S$-category $L^{H}(C,S)$, called the \textit{hammock localization}, which is an enhanced version of the localized category $W^{-1}C$. $L^{H}(C,W)$ (often
denoted simply by $L^{H}(C)$ when $W$ is clear from the context) is a model for the Dwyer-Kan simplicial
localization of $C$ along $W$ (\cite{dk1} and \cite{dk3}). $L^{H}(C,W)$ has the advantage, with respect to the simplicial localization of \cite{dk1}, that there is a natural morphism (called localization morphism) of $S$-categories $L:C\rightarrow L^{H}(C,W)$. With the same notations, we will write $\mathrm{Ho}(C)$ for the standard localization $W^{-1}C$ and call it the \textit{homotopy category} of $C$. In such a context, we will say that two objects
in $C$ are \emph{equivalent} if they are linked by a string of morphisms
in $W$. Equivalent objects in $C$ goes to isomorphic objects in $\mathrm{Ho}(C)$, but the
contrary is not correct in general (though it will be true in most of our
context, e.g. when $C$ is a model category or a good Waldhausen category).\\
The construction $(C,W) \mapsto L^{H}(C,W)$ is functorial in the pair $(C,W)$ and it also
extends naturally to the case where $W$ is a sub-$S$-category
of an $S$-category $C$ (see \cite[\S $6$]{dk1}).
Two fundamental properties of the functor $L^{H}:(C,W) \mapsto L^{H}(C,S)$ are the following
\begin{itemize}

\item The localization morphism $L$ identifies $\pi_{0}(L^{H}(C,W))$ with the (usual, Gabriel-Zisman) localization $W^{-1}C$.

\item If $M$ be a simplicial model category, and $W\subset M$ is its subcategory of equivalences, then
the full sub-$S$-category $M^{cf}$ of $M$, consisting of objects which are cofibrant and fibrant, is equivalent to $L^{H}(M,W)$.
\end{itemize}

We will neglect all kind of considerations about universes in our set-theoretic and categorical setup,
leaving to the reader to keep track
of the various choices of universes one needs in order the different constructions to make sense. \\

\end{section}

\begin{section}{Good Waldhausen categories}

In this section we introduce the class of Waldhausen categories (\textit{good} Waldhausen categories) we are going to work with and for
which our main theorem (Theorem \ref{t1}) holds. Regarding the choice of this class, it turns out in practice that, though some usual
Waldhausen categories might not be \textit{good} in our sense, to our knowledge there always exists a \textit{good Waldhausen
model} for them, i.e. a good Waldhausen category with the same $K$-theory space (up to homotopy). In other words, we do not know any relevant example
which, for $K$-theoretical purposes, could require using non-good Waldhausen categories.

We would also like to stress that the class of \textit{good} Waldhausen categories is not the most
general one for which our results hold. As the
reader will notice, our main results should still be correct for any Waldhausen category having a good
enough \textit{homotopy calculus of fractions} (in the sense of \cite[\S 6]{dk3}). \\

If $M$ is a model category, we denote by $M^{f}$ its full subcategory of fibrant objects. When the model
category $M$ is pointed,
the category $M^{f}$ will be considered as a Waldhausen category
in which weak equivalences and fibrations are induced by the model structure of $M$. \\

\begin{df}\label{d1}
A \emph{good Waldhausen category} is a Waldhausen category $C$ for which there
exists a pointed model category $M$ and a fully faithful functor $i : C \longrightarrow M^{f}$
satisfying the following conditions
\begin{enumerate}
\item The functor $i$ commutes with finite limits (in particular $i(*)=*$).
\item The essential image $i(C) \subset M^{f}$ is stable by weak equivalences (i.e., if $\;x \in M^{f}$ is
weakly equivalent to an object of $i(C)$, then $x \in i(C)$).
\item A morphism in $C$ is a fibration (resp. a weak equivalence) if and only if its image
is a fibration (resp. a weak equivalence) in $M$.
\end{enumerate}
\end{df}

\medskip

Most of the time we will identify $C$ with its essential image $i(C)$ in $M$ and forget about the
functor $i$. However, the model category $M$ and the embedding $i$ are \textit{not} part of the data. \\

\begin{ex}\label{nongood}
\emph{Let $k$ be a ring and $\mathrm{Ch}(k)$ be the category of (unbounded) chain complexes of $k$-modules.
The category $\mathrm{Ch}(k)$ is a model category with weak equivalences (resp. fibrations) given by the quasi-isomorphisms
(resp. by the epimorphisms). The subcategory $V$ of bounded complexes of finitely generated projective $k$-modules
is a Waldhausen category, where fibrations and weak equivalences are induced by $\mathrm{Ch}(k)$. However, $V$ might not be
a good Waldhausen category because it is not closed under quasi-isomorphisms in $\mathrm{Ch}(k)$: its closure
is the category $\mathrm{Perf}(k)$ of perfect complexes in $\mathrm{Ch}(k)$, which is indeed a good Waldhausen category
(for the induced structure). Nevertheless, the $K$-theory spectra of $V$ and of $\mathrm{Perf}(k)$ are
naturally equivalent. This is a typical situation of a Waldhausen category that might not be good
but which admits a good Waldhausen model.}
\end{ex}

It is clear from Definition \ref{d1} that
any morphism $f : x \rightarrow y$ in a good Waldhausen category $C$
possesses a (functorial) factorization
$$f : \xymatrix{x \ar[r]^-{j} & x' \ar[r]^-{p} & y,}$$
where $j$ is a cofibration and $p$ a fibration, and one of them is
a weak equivalence. Here, by cofibration in $C$ we mean a morphism
that has the left lifting property with respect to all fibrations
in $C$ that are also weak equivalences. Using
\cite[$8.2$]{dk3} (with $W_{1}$ being the class of
trivial cofibrations in $C$, and $W_{2}$ the class of trivial fibrations),
one sees that the existence of such factorizations
implies that the category $C$ has a two sided homotopy calculus of fractions
with respect to the weak equivalences $W$. In particular, the simplicial sets
of morphisms in $L^{H}C$ can be computed using hammocks of types
$W^{-1} C W^{-1}$ (\cite[Prop. $6.2$ (i)]{dk3}). As an immediate consequence, we get the following
important fact.

\begin{prop}\label{p1'}
Let $C$ be a good Waldhausen category, and $W$ its sub-category of
weak equivalences. Let $wL^{H}C$ be the sub-$S$-category
of $L^{H}C$ consisting of homotopy equivalences (i.e. of morphisms
projecting to isomorphisms in $\pi_{0}(L^{H}C)$). Then, the natural morphism
induced on the geometric realizations
$$|W| \longrightarrow |wL^{H}C|$$
is a weak equivalence of simplicial sets.
\end{prop}

\textit{Proof:} Indeed, as $C$ has a two sided homotopy calculus of fractions,
then \cite[$6.2$ $(i)$]{dk3} implies that the natural
morphism $|L^{H}W| \longrightarrow  |wL^{H}C|$ is a weak equivalence.
As
\cite[4.2]{dk1} and \cite[Prop. 2.2]{dk3} implies that the natural morphism
$|W| \longrightarrow |L^{H}W|$ is also a weak equivalence, so
is the composition
$$|W| \longrightarrow |L^{H}W| \longrightarrow |wL^{H}C|$$
\hfill $\Box$ \\

Let $C$ be a good Waldhausen category, and $i : C \hookrightarrow M^{f}$
an embedding as in Definition \ref{d1}.
Condition $(2)$ of Definition \ref{d1}, and the definition of the hammock localization
of \cite[2.1]{dk3}, implies immediately that the induced morphism of $S$-categories
$$L^{H}C \longrightarrow L^{H}M^{f}$$
is fully faithful (in the sense of defintion \ref{eqS}).
This implies that the (homotopy type of the) simplicial sets of morphisms of $L^{H}C$ can actually be computed
in the model category $M$ by using the standard simplicial and co-simplicial resolutions techniques available in model categories (see \cite{dk4}).

The induced functor $\mathrm{Ho}(C) \longrightarrow \mathrm{Ho}(M)$ being fully faithful, one sees that
any morphism $a\rightarrow b$ in the homotopy category of
a good Waldhausen category $C$ can be represented by a diagram $\xymatrix{a & \ar[l]_-{u} a' \ar[r] & b}$ in $C$,
where $u$ is a weak equivalence (recall that any object in $C$ is fibrant in $M$). From a general point of view, homotopy categories
of good Waldhausen categories behave very much like
categories of fibrant objects in a model categories. For example the set of morphisms in the
homotopy category can be computed using homotopy classes of morphisms
from cofibrant to fibrant objects (as explained in \cite{ho}).
In this work we will often use implicitely all these properties.

\begin{ex}\label{ex}
\emph{
\begin{enumerate}
\item
The first standard example of a good Waldhausen category is the category $\mathrm{Perf}(k)$ of
perfect complexes over a ring $k$. Recall that the fibrations are the epimorphisms and the quasi-isomorphisms
are the weak equivalences. The category $\mathrm{Perf}(k)$ is clearly a full subcategory of $\mathrm{Ch}(k)$, the
category of all chain complexes of $k$-modules. If we endow $\mathrm{Ch}(k)$ with its projective model structure
of \cite[Thm 2.3.11]{ho} (for which the weak equivalences are the quasi-isomorphisms and the fibrations are the epimorphisms)
then one checks immediately that the conditions of the Definition \ref{d1} are satisfied.
\item The previous example can be generalized in order to construct a good Waldhausen category
that computes the $K$-theory of schemes in the sense of \cite{tt}. One possible way to do this, is
by using the model category $\mathrm{Ch_{QCoh}}(X)$ of complexes of quasi-coherent $\mathcal{O}_{X}$-Modules on a
quasi-compact and quasi-separated scheme $X$ defined
in \cite[Cor. 2.3 (b)]{ho2}. Recall that in this injective model structure the cofibrations are the monomorphisms and weak equivalences are the
quasi-isomorphisms. Inside $\mathrm{Ch_{QCoh}}(X)$ we have the full subcategory of perfect complexes
$\mathrm{Perf}(X) \subset \mathrm{Ch_{QCoh}}(X)$, which is a Waldhausen category for which
weak equivalences are the quasi-isomorphisms and the fibrations are the epimorphisms, and
that computes the $K$-theory of the scheme $X$. This Waldhausen category does not seem
to be good in the sense of Definition \ref{d1}, however its full subcategory of fibrant
objects $\mathrm{Perf}(X)^{f} \subset \mathrm{Perf}(X)$ is good if we endow it with the induced structure of Waldhausen category
coming from $\mathrm{Ch_{QCoh}}(X)^{f}$. Now, the inclusion functor $\mathrm{Perf}(X)^{f} \hookrightarrow \mathrm{Perf}(X)$
is an exact functor of Waldhausen categories (as
fibrations in $\mathrm{Ch_{QCoh}}(X)$ are in particular epimorphisms, \cite[Prop. 2.12]{ho2}), and the approximation theorem (\cite[Theorem 1.6.7]{wa})
tells us that it induces a weak equivalence on the corresponding $K$-theory spectra. Therefore, the $K$-theory
of the scheme $X$ can be computed using the good Waldhausen category $\mathrm{Perf}(X)^{f}$.
\item More generally, for any ringed site $(C,\mathcal{O})$, there exists a good Waldhausen
category that computes the $K$-theory of the Waldhausen category of perfect complexes
of $\mathcal{O}$-modules on $C$. This requires a model category structure on the
category of complexes of $\mathcal{O}$-modules that we will not describe in this work.
\item For a topological space $X$, one can use the model category of spaces
under-and-over $X$ (retractive spaces over $X$), $X/\mathrm{Top}/X$, in order to define a good Waldhausen category
computing the Waldhausen $K$-theory of $X$ (\cite[2.1]{wa}).
\end{enumerate}
}

\end{ex}
\end{section}

\begin{section}{$DK$-equivalences and the approximation theorem}

If $C$ is any Waldhausen category we will simply denote
by $L^{H}C$ its hammock localization along the subcategory of weak equivalences
as defined in \cite[2.1]{dk3}. $L^{H}C$ is an $S$-category that comes together with a
localization functor
$$l : C \longrightarrow L^{H}C.$$
If $f : C \longrightarrow D$ is an exact functor between Waldhausen categories, then
it induces a well defined morphism of $S$-categories $Lf : L^{H}C \longrightarrow L^{H}D$,  such that the following diagram is commutative
$$\xymatrix{
C \ar[r]^-{l} \ar[d]_-{f} & L^{H}C \ar[d]^-{Lf} \\
D \ar[r]^-{l} & L^{H}D.}
$$

\begin{df}\label{d2}
An exact functor $f : C \longrightarrow D$ between Waldhausen categories
is a \emph{DK-equivalence} if the induced morphism
$Lf : L^{H}C \longrightarrow L^{H}D$ is an equivalence of
$S$-categories (in the sense of \emph{\cite[1.3 (ii)]{dk2}}).
\end{df}

\medskip

\noindent Obviously, the expression \textit{DK-equivalence} refers to Dwyer and Kan. \\

The following proposition is a strong form of the approximation
theorem for good Waldhausen categories. It is probably false
for more general Waldhausen categories.

\begin{prop}\label{p1}
If $f : C \longrightarrow D$ is a DK-equivalence between good Waldhausen categories
then the induced morphism on the $K$-theory spectra
$$K(f) : K(C) \longrightarrow K(D)$$
is a weak equivalence.
\end{prop}

\textit {Proof:} Let $S_{n}C$ and $S_{n}D$ denote the dual versions (with cofibration replaced by fibrations)
of the categories with weak equivalences defined and denoted in the same way in \cite[1.3]{wa}. We will prove the following more precise claim \\

\noindent \textsl{Claim:} ``For any $n\geq 0$, the induced functor
$$S_{n}f : S_{n}C \longrightarrow S_{n}D$$
induces a weak equivalence on the classifying spaces of weak equivalences
$$wS_{n}f : |wS_{n}C| \simeq |wS_{n}D|.\textrm{''}$$

\medskip

Note that
the category $wS_{n}C$ is equivalent to the category of strings of fibrations in $C$
$$\xymatrix{x_{n} \ar[r] & x_{n-1} \ar[r] & \dots \ar[r] & x_{1}}$$
and levelwise weak equivalences between them.
As nerves of categories are preserved (up to a weak equivalence) by equivalences of categories, we can
assume that $S_{n}C$ (resp. $S_{n}D$) \textit{actually is} the category of strings of fibrations in $C$
(resp., in $D$); the fact that $S_{n}C$ is a bit more complicated than
just the category of strings of fibrations is only used to have
a strict simplicial diagram of categories $[n] \mapsto S_{ n}C$ (see \cite[1.3, p. 329]{wa}).

\begin{lem}\label{l1}
Let $f : C \longrightarrow D$ be a DK-equivalence between good Waldhausen categories.
Then the induced morphism on the classifying spaces
$$|f| : |wC| \longrightarrow |wD|$$
is a weak equivalence.
\end{lem}

\textit{Proof of Lemma.} This follows immediately from Proposition \ref{p1'}. \hfill $\Box$ \\

Note that the previous lemma already implies the \textsl{Claim} above for $n=1$. For general $n$, it is
then enough to prove that the categories $S_{n}C$ and
$S_{n}D$ are again good Waldhausen categories and that
the induced exact functor
$$S_{n}f : S_{n}C \longrightarrow S_{n}D$$
is again a DK-equivalence, and then apply Lemma \ref{l1} to get the \textsl{Claim}.

We need to recall here the Waldhausen structure on the category
$S_{n}C$. The fibrations (resp. weak equivalences) are the morphisms
$$\xymatrix{
x_{n} \ar[r] \ar[d] & x_{n-1} \ar[r] \ar[d] & \dots \ar[r] & x_{1} \ar[d] \\
y_{n} \ar[r]  & y_{n-1} \ar[r] & \dots \ar[r] & y_{1}}$$
such that each induced morphism
$$x_{i} \longrightarrow x_{i-1}\times_{y_{i-1}}y_{i}$$
is a fibration in $C$ (resp., such that each morphism
$x_{i} \rightarrow y_{i}$ is a weak equivalence in $C$).
With this definition we have

\begin{lem}\label{l2}
If $C$ is a good Waldhausen category then so is $S_{n}C$.
\end{lem}

\textit{Proof of lemma.}
Let us consider an embedding $C \subset M^{f}$, of $C$
in the category of fibrant objects in a pointed model category $M$ (and
satisfying the conditions of Definition \ref{d1}). We consider the
category $S_{n}M:=M^{I(n-1)}$, of strings of $(n-1)$ composable morphisms in $M$.
Here we have denoted by $I(n-1)$ the free category with $n$ composable morphism
$$I(n-1):=\left\{ \xymatrix{n \ar[r] & n-1 \ar[r] & \dots \ar[r] & 1} \right\}.$$

The objects of $S_{n}M$ are therefore diagrams in $M$
$$\xymatrix{x_{n} \ar[r] & x_{n-1} \ar[r] & \dots \ar[r] & x_{1}.}$$
We endow the category $S_{n}M$ with the model structure for which weak equivalences
are defined levelwise. The fibrations are morphisms
$$\xymatrix{
x_{n} \ar[r] \ar[d] & x_{n-1} \ar[r] \ar[d] & \dots \ar[r] & x_{1} \ar[d] \\
y_{n} \ar[r]  & y_{n-1} \ar[r] & \dots \ar[r] & y_{1}}$$
such that each induced morphism
$$x_{i} \longrightarrow x_{i-1}\times_{y_{i-1}}y_{i}$$
is a fibration in $M$. Note that in particular
fibrant objects in $S_{n}M$ are strings of fibrations in $M$
$$\xymatrix{x_{n} \ar[r] & x_{n-1} \ar[r] & \dots \ar[r] & x_{1}.}$$
This model structure is known as the Reedy model structure
described e.g. in \cite[Thm. 5.2.5]{ho}, when the category $I(n-1)$ is considered as a
Reedy category in the obvious way.
Now, the category $S_{n}C$ has an induced natural embedding into
$(S_{n}M)^{f}$, which satisfies the conditions of Definition \ref{d1}.
This concludes the proof of Lemma \ref{l2}. \hfill $\Box$

\begin{lem}\label{l2'}
The induced exact functor
$$S_{n}f : S_{n}C \longrightarrow S_{n}D$$
is a DK-equivalence.
\end{lem}

\textit{Proof of lemma.} Let us first show that the induced morphism
$$L^{H}S_{n}f : L^{H}S_{n}C \longrightarrow L^{H}S_{n}D$$
is fully faithful. To see this, let $C \hookrightarrow M^{f}$ be
an embedding of $C$ in a pointed model category as in Definition \ref{d1}.
Then, the simplicial sets of morphisms in $L^{H}C$ are equivalent to
the corresponding mapping spaces computed in the model category $M$. Applying this
argument to the embedding $S_{n}C \hookrightarrow S_{n}M$, we deduce that the
simplicial sets of morphisms of $L^{H}S_{n}C$ are equivalent to the corresponding mapping spaces
computed in the model category $S_{n}M$. Finally, it is quite easy to
compute the mapping spaces in $S_{n}M$ in terms of the mapping spaces
of $M$.  The reader will
check that
the simplicial set of morphisms from $\xymatrix{x:=(x_{n} \ar[r] & x_{n-1} \ar[r] & \dots \ar[r] & x_{1})}$
to $\xymatrix{y:= (y_{n} \ar[r] & y_{n-1} \ar[r] & \dots \ar[r] & y_{1})}$ in
$L^{H}S_{n}C$ is given by the following iterated homotopy fiber product
$$ \underline{Hom}_{L^{H}C}(x_{n},y_{n})\times^{h}_{\underline{Hom}_{L^{H}C}(x_{n},y_{n-1})}
\underline{Hom}_{L^{H}C}(x_{n-1},y_{n-1})\times^{h} \dots
\times^{h}_{\underline{Hom}_{L^{H}C}(x_{2},y_{1})}\underline{Hom}_{L^{H}C}(x_{1},y_{1}).$$
This description of the simplicial sets of morphisms in $L^{H}S_{n}C$ is of course
also valid for $L^{H}S_{n}D$. It shows in particular that if $L^{H}f : L^{H}C \longrightarrow
L^{H}D$ is fully faithful, then so is
$L^{H}S_{n}f : L^{H}S_{n}C \longrightarrow L^{H}S_{n}D$ for any $n$. \\

It remains to show that the morphism
$$L^{H}S_{n}f : L^{H}S_{n}C \longrightarrow L^{H}S_{n}D$$
is essentially surjective. It is enough
to prove that for any object
$$\xymatrix{y:= (y_{n} \ar[r] & y_{n-1} \ar[r] & \dots \ar[r] & y_{1})}$$
\noindent in $S_{n}D$, there exists an object $$\xymatrix{x:= (x_{n} \ar[r] & x_{n-1} \ar[r] & \dots \ar[r] & x_{1})}$$
\noindent in $S_{n}C$, an object $$\xymatrix{z:= (z_{n} \ar[r] & z_{n-1} \ar[r] & \dots \ar[r] & z_{1})}$$
\noindent in $S_{n}D$ and a diagram of weak equivalences in $S_{n}D$
$$\xymatrix{f(x) & \ar[l] z \ar[r] & y.}$$
For this, we let $z \rightarrow y$ be a cofibrant replacement of $y$ in the good
Waldhausen category $S_{n}D$ (recall that cofibrations in
a Waldhausen category are defined to be morphisms having the left lifting property with
respect to fibrations which are weak equivalences; by definition of a good Waldhausen category,
a cofibrant replacement functor always exists).
By induction, we may assume that there
exists $\xymatrix{x_{\leq n-1}:= (x_{n-1}\ar[r] & x_{n-2} \ar[r]  & \dots \ar[r] & x_{1})} \in S_{n-1}D$, and
a weak equivalence $z_{\leq n-1} \rightarrow f(x_{\leq n-1})$ in $S_{n-1}D$,
where $\xymatrix{z_{\leq n-1}:= (z_{n-1}\ar[r] & z_{n-2} \ar[r]  & \dots \ar[r] & z_{1})} \in S_{n-1}D$.
And it remains to show that there exists a fibration
$x_{n} \rightarrow x_{n-1}$ in $C$, and a weak equivalence $z_{n} \rightarrow f(x_{n})$ in $D$, such that the
following diagram in $D$ commutes
$$\xymatrix{
z_{n} \ar[r] \ar[d] & z_{n-1} \ar[d] \\
f(x_{n}) \ar[r] & f(x_{n-1}).}$$
As $f$ induces an equivalence $\mathrm{Ho}(C)\simeq \mathrm{Ho}(D)$ on the homotopy categories and the $z_{i}$ are cofibrant objects in $D$,
it is clear that one can find a fibration $x_{n} \rightarrow x_{n-1}$ in $C$, and a weak equivalence $z_{n} \rightarrow f(x_{n})$ in $D$
such that the above diagram is commutative in $\mathrm{Ho}(D)$. But, as $z_{n}$ is cofibrant
and $f(x_{n}) \rightarrow f(x_{n-1})$ is a fibration between fibrant objects, we can always
choose the weak equivalence $z_{n} \rightarrow f(x_{n})$ in such a way that the above diagram
commutes in $D$ (the argument is the same as in the case of model categories, therefore we leave the details to the reader).

This construction gives the required diagram in $D$
$$\xymatrix{f(x) & \ar[l] z \ar[r] & y,}$$
and concludes the proof of the lemma. \hfill $\Box$ \\

Lemmas \ref{l2} and \ref{l2'} show that $S_{n}f : S_{n}C \longrightarrow S_{n}D$ is
also a DK-equivalence between good Waldhausen categories
for any $n$, and
therefore Lemma \ref{l1} finishes the proof of the \textsl{Claim} and therefore of Proposition \ref{p1}.
\hfill $\Box$

\end{section}

\begin{section}{Simplicial localization of good Waldhausen categories}

Given fibrant simplicial sets $X$, $Y$ and $Z$, and a diagram
$\xymatrix{X \ar[r] & Z & Y \ar[l]}$, we denote by
$X\times_{Z}^{h}Y$ the corresponding standard homotopy fibered product. Explicitly, it is
defined by
$$X\times_{Z}^{h}Y:=(X\times Y)\times_{Z\times Z}Z^{\Delta^{1}}.$$
Note that for any simplicial set $A$, there is a natural isomorphism of
simplicial sets
$$\underline{Hom}(A,X\times_{Z}^{h}Y)\simeq
\underline{Hom}(A,X)\times_{\underline{Hom}(A,Z)}^{h}\underline{Hom}(A,Y).$$

\begin{df}\label{d3}
Let $T$ be an $S$-category.
We say that $T$ is \emph{pointed} if there exists an object $* \in T$
such that for any other object $x\in T$, the simplicial sets
$\underline{Hom}_{T}(x,*)$ and $\underline{Hom}_{T}(*,x)$
are weakly equivalent to $*$.
\end{df}

For the next definition, recall that an $S$-category is said to be \textit{fibrant} if all
its simplicial sets of morphisms are fibrant simplicial sets. The existence of a model
structure on $S$-categories with a fixed set of objects (see for example \cite{dk1}) implies
that for any $S$-category $T$, there exists a fibrant $S$-category $T'$ and an
equivalence of $S$-categories $T \longrightarrow T'$ (this equivalence is furthermore the identity on the
set of objects). Such a $T'$ will be called a \textit{fibrant model of $T$}.

\begin{df}\label{d4}
\begin{enumerate}
\item Let $T$ be a fibrant $S$-category. We say that $T$ \emph{has fibered products} if for any diagram of
morphisms in $T$
$$\xymatrix{x \ar[r]^-{u} & z & \ar[l]_-{v} y}$$
there exists an object $t \in T$, two morphisms
$$\xymatrix{x & \ar[l]_-{p} t \ar[r]^-{q} & y}$$
and a homotopy $h \in \underline{Hom}_{T}(t,z)^{\Delta^{1}}$
such that
$$\partial_{0}h=u\circ p \qquad \partial_{1}h=v\circ q,$$
and which satisfies the following
universal property:

``for any object $w \in T$, the natural morphism induced by $(p,q,h)$
$$\underline{Hom}_{T}(w,t) \longrightarrow \underline{Hom}_{T}(w,x)\times^{h}_{\underline{Hom}_{T}(w,z)}
\underline{Hom}_{T}(w,y)$$
is a weak equivalence.''

Such an object $t$ together with the data $(p,q,h)$ is called
a \emph{fibered product} of the diagram $\xymatrix{x\ar[r] & z & \ar[l] y}$.

\item
For a general $S$-category $T$, we say that $T$ \emph{has fibered products} if
one of its fibrant model has.
\end{enumerate}
\end{df}

The fact that most of the time the hammock localization of a category
is not a fibrant $S$-category might be annoying. There exists however a
more intrinsic version of the above definition that we now describe.

For an $S$-category $T$, one can consider
the category of simplicial functors from its opposite $S$-category $T^{op}$ to the
category $SSet$, of simplicial sets. This category is a model category for which
the weak equivalences and the fibrations are defined levelwise; we will denote it by
$\mathrm{SPr}(T)$. There exists a simplicially-enriched Yoneda functor
$$\underline{h} : T \longrightarrow \mathrm{SPr}(T)$$
that sends an object $x \in T$ to the diagram
$$\begin{array}{cccc}
\underline{h}_{x} : & T^{op} & \longrightarrow & SSet \\
 & y & \mapsto & \underline{Hom}_{T}(y,x).
\end{array}$$
We will say that an object of $\mathrm{SPr}(T)$ is \textit{representable} if it is weakly equivalent in $\mathrm{SPr}(T)$ to
some $\underline{h}_{x}$ for some $x\in T$. With these notions
the reader will check the following fact as an exercise.

\begin{lem}\label{l3}
An $S$-category $T$ has fibered products if and only if the
full subcategory of $\mathrm{SPr}(T)$ consisting of representable objects is stable under
homotopy pull-backs.
\end{lem}

Note that one can then assume the property after the ``iff'' in the previous Lemma as an  equivalent
definition of $S$-category with fiber products. The following proposition is well known when $C$ is a model category (see for example \cite[8.4]{sh}).

\begin{prop}\label{p2}
Let $C$ be a good Waldhausen category. Then the $S$-category $L^{H}C$
is pointed and has fibered products.
\end{prop}

\textit {Proof:} Let $C \hookrightarrow M$ be an embedding of $C$ as
a full subcategory of a pointed model category as in Definition \ref{d1}. The conditions
of \ref{d1} imply that the induced morphism of $S$-categories
$$L^{H}C \longrightarrow L^{H}M$$
is fully faithful (i.e. induces a weak equivalence on the corresponding simplicial sets of morphisms).
Therefore, $L^{H}C$ is equivalent to the full sub-$S$-category of $L^{H}M$ consisting
of objects belonging to $C$. Now, it is well known that the $S$-category
$L^{H}M$ has fibered products and furthermore that the fibered products in $L^{H}M$
can be identified with the homotopy fibered products in the model category $M$ (see \cite[8.4]{sh}).
By Definition \ref{d1}, the full sub-$S$-category of $L^{H}M$
of objects belonging to $C$ is therefore stable by fibered products. This formally implies that
$L^{H}C$ has fibered products (the details are left to the reader).

In the same way, as $M$ is a pointed model category, the object $*$ in $M$, viewed as an object in $L^{H}M$,
satisfies the condition of Definition \ref{d3}, so $L^{H}M$ is a pointed $S$-category.
But, by condition $(1)$ of Definition \ref{d1}, this $*$ belongs to the image
of $L^{H}C$ in $L^{H}M$. As $L^{H}C \longrightarrow L^{H}M$ is fully faithful, this
shows that $L^{H}C$ is a pointed $S$-category. \hfill $\Box$

\end{section}

\begin{section}{$K$-theory of $S$-categories}

In this section, we define for any pointed $S$-category $T$ with fibered products
a $K$-theory spectrum $K(T)$. We will show that
$K(T)$ is invariant, up to weak equivalences, under equivalences of $S$-categories in $T$. The construction $T \mapsto K(T)$ is also functorial
in $T$, but we will not investigate this in this work, as it is more technical to prove and is not really needed for our main purpose. \\

We fix $T$ a \textit{pointed} $S$-category \textit{with fibered products}. We consider
the model category $\mathrm{SPr}(T)$ of simplicial diagrams on $T^{op}$, and its
associated Yoneda embedding
$$\begin{array}{cccc}
\underline{h} : & T  &\longrightarrow & \mathrm{SPr}(T) \\
 & x & \mapsto & \underline{Hom}_{T}(-,x)
\end{array}$$

Recall the following homotopy version of the simplicially enriched Yoneda lemma (e.g. \cite[Prop. 2.4.2]{tv})

\begin{lem}\label{hSYoneda}
Let $T'$ be any $S$-category. For any object $F\in \mathrm{SPr}(T')$ and any object $x \in T'$, there
is a natural isomorphism in the homotopy catgory of simplicial sets
$$\mathbb{R}\underline{Hom}(\underline{h}_{x},F)\simeq F(x).$$
\noindent In particular, the induced functor $\underline{h}:\pi_{0}T' \rightarrow \mathrm{Ho}(\mathrm{SPr}(T'))$ is fully faithful.
\end{lem}

Recall that any object in $\mathrm{SPr}(T)$ which is weakly equivalent to
some $\underline{h}_{x}$ is called \textit{representable}.
If $*$ denotes the final object in $\mathrm{SPr}(T)$, let us consider the model category
$$\widehat{T}_{*}:=*/\mathrm{SPr}(T),$$
of pointed objects in the model category $\mathrm{SPr}(T)$ (see \cite[p. 4 and Ch. 6]{ho}).
Clearly, $\widehat{T}_{*}$ is a pointed model category. We now
consider its full subcategory of fibrant objects,
denoted by $\widehat{T}_{*}^{f}$, and define the category $M(T)$ to be the
full subcategory consisting of objects in $\widehat{T}^{f}_{*}$ whose underlying objects in
$\mathrm{SPr}(T)$ are representable.

As we supposed that $T$ has fibered products, one checks immediately that
$M(T)$ is a full subcategory of $\widehat{T}^{f}_{*}$ which is
stable under weak equivalences and homotopy pull-backs (see Lemma \ref{l3}). Moreover, $M(T)$ contains
the final object of $\widehat{T}_{*}$ since this is weakly equivalent to $\underline{h}_{*}$ for
any object $* \in T$ as in Definition \ref{d4}. Therefore, endowed with the induced Waldhausen structure
coming from $\widehat{T}_{*}^{f}$, $M(T)$ clearly becomes a good Waldhausan category.

\begin{df}\label{d5}
The \emph{$K$-theory spectrum of the $S$-category $T$} is defined to be the $K$-theory spectrum of
the Waldhausen category $M(T)$. It is denoted by
$$K(T):=K(M(T)).$$
\end{df}

We will now show that the construction $T \mapsto K(T)$ is
functorial \textit{with respect to equivalences} of $S$-categories. Though $T \mapsto K(T)$ actually satisfies a
more general functoriality property, its functoriality with respect to equivalences of $S$-categories will be
enough for our present purpose which is to deduce that $K(T)$ only depends (up to weak equivalences) on the $S$-equivalence class of $T$.\\

Let $f : T \longrightarrow T'$ be an equivalence of pointed $S$-categories with fibered products.
We deduce a pull-back functor
$$f^{*} : \mathrm{SPr}(T') \longrightarrow \mathrm{SPr}(T),$$
as well as its pointed version
$$f^{*} : \widehat{T'}_{*} \longrightarrow \widehat{T}_{*}.$$
This functor is in fact a right Quillen functor whose letf adjoint is denoted
by
$$f_{!} : \widehat{T}_{*} \longrightarrow \widehat{T'}_{*}.$$
As the morphism $f$ is an equivalence of $S$-categories, this Quillen adjunction
is actually known to be a Quillen equivalence (see \cite{dk2}). The functor $f^{*}$ (pointed-version) being right Quillen, it induces
a functor on the subcategories of fibrant objects
$$f^{*} : \widehat{T'}^{f}_{*} \longrightarrow \widehat{T}^{f}_{*}.$$

\begin{prop}\label{p4}
The functor above sends the subcategory $M(T') \subset (\widehat{T'}_{*})^{f}$
into the subcategory $M(T) \subset (\widehat{T}_{*})^{f}$.
\end{prop}

\textit {Proof:} By definition of $M(-)$, it is enough to show that the right derived functor
$$\mathbb{R}f^{*}\simeq f^{*} : \mathrm{Ho}(\mathrm{SPr}(T')) \longrightarrow \mathrm{Ho}(\mathrm{SPr}(T))$$
preserves the property of being a representable object. But, this functor is an equivalence
of categories whose inverse is the functor
$$\mathbb{L}f_{!} : \mathrm{Ho}(\mathrm{SPr}(T)) \longrightarrow \mathrm{Ho}(\mathrm{SPr}(T')).$$
The reader will check that, by adjunction, one has for any object
$x \in T$ a natural isomorphism in $\mathrm{Ho}(\mathrm{SPr}(T'))$
$$\mathbb{L}f_{!}(\underline{h}_{x})\simeq \underline{h}_{f(x)}.$$
As $f$ is an equivalence of $S$-categories, for any object $y \in T'$ there exists
$x \in T$ and a morphism $u : f(x) \longrightarrow y$ in $T'$ inducing an isomorphism in $\pi_{0}T'$.
Clearly $\underline{h}_{u} : \underline{h}_{f(x)} \longrightarrow \underline{h}_{y}$
is an equivalence in $\mathrm{SPr}(T)$.
Therefore, one has
$$f^{*}(\underline{h}_{y})\simeq
f^{*}(\underline{h}_{f(x)})\simeq
f^{*}\circ \mathbb{L}f_{!}(\underline{h}_{x})\simeq \underline{h}_{x}.$$
This implies that $f^{*}\simeq \mathbb{R}f^{*}$ preserves representable objects. \hfill $\Box$ \\

The above proposition implies that any equivalence of $S$-categories $f:T \rightarrow T'$ induces a well defined exact functor of good Waldhausen categories
$$f^{*} : M(T') \longrightarrow M(T).$$
The rule $f \mapsto f^{*}$ is clearly controvariantly functorial in $f$ (i.e. one has a natural isomorphism
$(g\circ f)^{*}\simeq f^{*}\circ g^{*}$, satisfying the usual co-cycle condition).
Therefore we get a controvariant (lax) functor from the category of pointed $S$-categories with fibered products and
$S$-equivalences to the category of Waldhausen categories and exact functors\footnote{By applying the
standard strictification procedure we will assume that $M \mapsto M(T)$ is a genuine functor from
$S$-categories towards Waldhausen categories}.

\begin{prop}\label{p5}
Let $f : T \longrightarrow T'$ be an equivalence of pointed $S$-categories with
fibered products. Then the induced exact functor
$$f^{*} : M(T') \longrightarrow M(T)$$
is a DK-equivalence (see Definition \ref{d2}).
\end{prop}

\textit {Proof:} By construction, there is a commutative diagram on the level
of hammock localizations
$$\xymatrix{
L^{H}(\widehat{T'}_{*}^{f}) \ar[r]^-{L^{H}f^{*}} & L^{H}(\widehat{T}^{f}_{*}) \\
L^{H}M(T') \ar[r] \ar[u] & L^{H}M(T) \ar[u].}$$
The functor $f^{*}$ being a Quillen equivalence it is well known that
the top horizontal arrow is an equivalence of $S$-categories (\cite{dk3}). But, as the
vertical morphisms of $S$-categories are fully faithful this implies that the morphism
$$L^{H}f^{*} : L^{H}M(T') \longrightarrow L^{H}M(T)$$
is fully faithful. But the isomorphism in $\mathrm{Ho}(\mathrm{SPr}(T))$
$$f^{*}(\underline{h}_{f(x)})\simeq \underline{h}_{x}$$
shows that the induced functor
$$\pi_{0}L^{H}f^{*} : \pi_{0}L^{H}M(T') \longrightarrow \pi_{0}L^{H}M(T)$$
is also essentially surjective,
and we conclude. \hfill $\Box$ \\

Using propositions \ref{p1}, \ref{p4} and \ref{p5}, we obtain the following conclusion.
Let us denote by $S-\mathrm{Cat}_{*}^{ex}$ the category of $S$-categories which are pointed and
have fibered products. Restricting the morphisms to equivalences of $S$-categories, we get  a subcategory
$wS-\mathrm{Cat}_{*}^{ex}$. Moreover, we denote by $\mathrm{Sp}$ the category of spectra, and
by $w\mathrm{Sp}$ its subcategory of weak equivalences. The previous constructions yield a well defined functor
$$\begin{array}{cccc}
K : & wS-\mathrm{Cat}_{*}^{ex} & \longrightarrow & w\mathrm{Sp}^{op} \\
& T & \mapsto & K(T):=K(M(T)) \\
& f & \mapsto & K(f^{*})
\end{array}$$

We can geometrically realize this functor to get a morphism on the corresponding classifying spaces
$$K : |wS-\mathrm{Cat}_{*}^{ex}| \longrightarrow |w\mathrm{Sp}^{op}|\simeq |w\mathrm{Sp}|,$$
which has to be understood as our $K$-theory functor from the
moduli space of pointed $S$-categories with fibered products to the moduli
space of spectra.

The fundamental groupoids of the spaces $|wS-\mathrm{Cat}_{*}^{ex}|$ and $|w\mathrm{Sp}|$
have the following description. Let us denote by $\mathrm{Ho}(S-\mathrm{Cat})$ (resp. by
$\mathrm{Ho}(\mathrm{Sp})$) the homotopy category of $S$-categories obtained by formally inverting
the $S$-equivalences (resp., the homotopy category
of spectra). Then, the fundamental groupoid $\Pi_{1}(|wS-\mathrm{Cat}_{*}^{ex}|)$ is naturally equivalent to the
sub-groupoid of $\mathrm{Ho}(S-\mathrm{Cat})$ consisting of pointed $S$-categories with fibered products and
isomorphisms between them (in $\mathrm{Ho}(S-\mathrm{Cat})$). In the same way, the fundamental groupoid
$\Pi_{1}(|w\mathrm{Sp}|)$ is naturally equivalent to the maximal sub-groupoid of $\mathrm{Ho}(\mathrm{Sp})$ consisting
of spectra and isomorphisms (in $\mathrm{Ho}(\mathrm{Sp})$). The $K$-theory morphism
$$K : |wS-\mathrm{Cat}_{*}^{ex}| \longrightarrow |w\mathrm{Sp}^{op}|\simeq |w\mathrm{Sp}|$$
\noindent defined above, induces
a well defined functor between the corresponding fundamental groupoids
$$K : \Pi_{1}(|wS-\mathrm{Cat}_{*}^{ex}|) \longrightarrow
\Pi_{1}(|w\mathrm{Sp}|)^{op}.$$
In other words, for any pair of pointed $S$-categories with fibered products $T$ and $T'$, and any
isomorphism $f : T \simeq T'$ in $\mathrm{Ho}(S-\mathrm{Cat})$, we have an isomorphism
$$K_{f}:=K(f^{*})^{-1} : K(T) \simeq K(T')$$
which is functorial in $f$.

Note however that the morphism $K : |wS-\mathrm{Cat}_{*}^{ex}| \longrightarrow |w\mathrm{Sp}|$
contains more information as for example it encodes the various morphisms on the simplicial monoids
of self-equivalences
$$\underline{Aut}(T) \longrightarrow \underline{Aut}(K(T)).$$

\begin{rmk} \emph{In closing this section, we would like to mention that the above construction of the $K$-theory spectrum $K(T)$ of an $S$-category $T$
can actually be made functorial enough in order to produce a well defined functor at the level of the underlying \textit{homotopy categories}
$$K : \mathrm{Ho}(S-\mathrm{Cat}_{*}^{ex}) \longrightarrow \mathrm{Ho}(\mathrm{Sp}).$$
Moreover, one can actually show that this can also be lifted to a \textit{morphism of} $S$\textit{-categories}
$$K : L^{H}(S-\mathrm{Cat}_{*}^{ex}) \longrightarrow L^{H}\mathrm{Sp}$$
between the corresponding hammock localizations, which is the best possible functoriality one could ever need in general.}
\end{rmk}

\end{section}

\begin{section}{Comparison}

In this section we prove that the $K$-theory spectrum of a good Waldhausen category (Definition \ref{d5}) can be
reconstructed from its simplicial localization. The main result is the following.

\begin{thm}\label{t1}
Let $C$ be a good Waldhausen category and $L^{H}C$ its hammock localization.
Then, there exists an isomorphism in the homotopy category of spectra
$$K(C) \simeq K(L^{H}C),$$
where the left hand side is the Waldhausen construction and the right hand side
is defined in Definition \ref{d5}.
\end{thm}

\textit {Proof:} We will explicitly produce a natural string of exact functors
between good Waldhausen categories, all of which are DK-equivalences, that links $C$ to $M(L^{H}C)$. Then, Proposition
\ref{p1} will imply the theorem. \\

We start by choosing a pointed model category $M$ and an embedding $C \hookrightarrow M^{f}$ as
in Definition \ref{d1}. Let $\Gamma$ be a cofibrant replacement functor in $M$, in the
sense of \cite[Def. 17.1.8]{hi}. Recall that this means that $\Gamma$ is a functor from $M$ to the category of
co-simplicial objects in $M$, together with a natural transformation $\Gamma \longrightarrow c$,
where $c$ is the constant co-simplicial diagram in $M$; moreover, for any
$x \in M$, the natural morphism
$$\Gamma(x) \longrightarrow c(x)$$
is a Reedy cofibrant replacement of the constant co-simplicial diagram $c(x)$
(i.e. it is a Reedy trivial fibration and $\Gamma(x)$ is cofibrant in the Reedy model
category (\cite[5.2]{ho}) of co-simplicial objects in $M$).

One should notice that if $x \in C$, since all the objects $\Gamma(x)^{n}$ are fibrant
objetcs in $M$ which are weakly equivalent to $x$, then $\Gamma(x)$ is actually
a co-simplicial object \textit{in} $C$.

Let us denote by $\widehat{C}$ the category of simplicial presheaves on $C$, and
by $\widehat{C}_{*}$ the category of pointed objects in $\widehat{C}$ (i.e. the category
of presheaves of pointed simplicial sets). Both these categories will be endowed with their
projective model structures for which fibrations and weak equivalences are defined objectwise.

For $x\in C$, we define a pointed simplicial presheaf
$$\begin{array}{cccc}
\underline{h}_{x} : & C^{op} & \longrightarrow & SSet_{*} \\
& y & \mapsto & \underline{h}_{x}(y):=Hom(\Gamma(y),x).
\end{array}$$
Note that $\underline{h}_{x}$ is a pointed simplicial presheaf because $C$ is pointed (and therefore the final object
 in $\widehat{C}$ can be identified with $\underline{h}_{*}$, where $*$ is the final \textit{and} the initial object in $C$).
The construction $x \mapsto \underline{h}_{x}$ then gives rise to a functor
$$\begin{array}{cccc}
\underline{h} : & C & \longrightarrow & \widehat{C}_{*} \\
& x & \mapsto & \underline{h}_{x}:=Hom(\Gamma(-),x).
\end{array}$$
As all objects in $C$ are fibrant in $M$, the standard properties
of mapping spaces tell us that for any $x\in C$ the pointed simplicial presheaf
$\underline{h}_{x}$ is a fibrant object in $\widehat{C}_{*}$ (see \cite[Cor. 17.5.3 (1)]{hi}). What we actually get is, therefore, a functor
$$\underline{h} : C \longrightarrow \widehat{C}_{*}^{f}.$$
If we endow the category $\widehat{C}_{*}^{f}$ with the induced Waldhausen structure coming from the projective model structure on $\widehat{C}_{*}$, the properties
of mapping spaces also imply that
the functor $\underline{h}$ is an exact functor between good Waldhausen categories (see \cite[Cor. 17.5.4 (2) and Cor. 17.5.5 (2)]{hi}).

We denote by $R(C)$ the full subcategory of $\widehat{C}_{*}^{f}$ consisting of objects weakly equivalent (in $\widehat{C}_{*}$) to
$\underline{h}_{x}$, for some $x\in C$. Objects in $R(C)$ will simply be called \textit{representable objects}.
As the functor $\underline{h}$ commutes with finite
limits, this subcategory is clearly a good Waldhausen category
when endowed with the induced Waldhausen structure.

\begin{lem}\label{l4}
The exact functor
between good Waldhausen categories $\underline{h} : C \longrightarrow R(C)$ is a DK-equivalence.
\end{lem}

\textit {Proof:} By construction the functor is essentially surjective up to weak equivalence, which implies
that $L^{H}\underline{h} : L^{H}C \longrightarrow L^{H}R(C)$ is indeed essentially surjective. It remains
to show that it is also fully faithful. Let us consider the composition
$$L^{H}C \longrightarrow L^{H}R(C) \longrightarrow L^{H}\widehat{C}_{*}^{f}.$$
The second morphism being fully faithful (as $R(C)$ is closed by weak equivalences in $\widehat{C}_{*}^{f}$), it is enough to show that
the composite morphism $L^{H}C \longrightarrow L^{H}\widehat{C}_{*}$ is fully faithful.
This essentially follows from the Yoneda lemma for pseudo-model categories of \cite[Lemma 4.2.2]{tv}, with the
small difference that $C$ is not exactly a pseudo-model category (\cite[Def. 4.1.1]{tv}), but only
the subcategory of fibrant objects in a pseudo-model category.

To fix this, we proceed as follows. Let $C'$ be the full subcategory of $M$ of objects weakly equivalent to some object in $C$.
Then, clearly $C'$ is a pseudo-model category (\cite[Def. 4.1.1]{tv}), there is an obvious embedding $C\hookrightarrow C'$
and (identifying $C$ with its essential image in $M$) its subcategory of fibrant objects $(C')^{f}$ coincides with $C$ .
Moreover, if we denote by $R$ a fibrant replacement functor in $M$, the functor $$\underline{h}_{R}:C' \longrightarrow
\widehat{C}_{*}$$ sending $x \in C'$ to $\underline{h}_{R(x)}$ preserves weak equivalences (\cite[Cor. 17.5.4 (2)]{hi}) and one has, by definition, a commutative diagram
$$\xymatrix{
C \ar[r]^-{\underline{h}} & \widehat{C}_{*} \\
C', \ar[ru]_-{\underline{h}_{R}} \ar[u]_{R} & }$$
giving rise to a corresponding commutative diagram of $S$-categories
$$\xymatrix{
L^{H}C \ar[r]^-{L^{H}\underline{h}}  & L^{H}\widehat{C}_{*} \\
L^{H}C'. \ar[ru]_-{L^{H}\underline{h}_{R}} \ar[u]^{L^{H}R} & }$$
Applying $L^{H}$ to the inclusion $C \hookrightarrow C'$, one sees that the morphism
$L^{H}R:L^{H}C' \longrightarrow L^{H}C$ is an equivalence of $S$-categories.
Finally, the
morphism $L^{H}\underline{h}_{R}:L^{H}C' \longrightarrow L^{H}\widehat{C}_{*}$ is fully faithful
by the following application of the Yoneda lemma for pseudo-model categories \cite[Lemma 4.2.2]{tv}.
For any $x$ and $y$ in $C'$, we have a chain of weak equivalences of simplicial sets
$$\underline{Hom}_{L^{H}(C')}(x,y)\simeq \mathrm{Map}_{M}(x,y)\simeq \mathrm{Map}_{M}(Rx,Ry) \simeq \underline{h}_{Ry}(Rx),$$
\noindent where $\mathrm{Map}(-,-)$ denotes the mapping space. But, by the standard  simplicially enriched Yoneda lemma, the
simplicial set $\underline{h}_{Ry}(Rx)$ is isomorphic to $\underline{Hom}_{\widehat{C'}}(h_{Rx},\underline{h}_{Ry})$, where
$h_{Rx}$ denotes presheaf of constant simplicial sets $z \mapsto Hom_{M}(z,Rx)$; moreover, if we let $W'$ denote the
weak equivalences in $C'$, $h_{Rx}$ is cofibrant $(C',W')^{\wedge}$ (defined in \cite[Def. 4.1.4]{tv}) and
$\underline{h}_{Ry}$ is fibrant in  $(C',W')^{\wedge}$. Hence $\underline{h}_{Ry}(Rx)$ is weakly equivalent to
$\mathrm{Map}_{(C',W')^{\wedge}}(h_{Rx},\underline{h}_{Ry})$ and then, by \cite[Lemma 4.2.2]{tv}, to
$\mathrm{Map}_{(C',W')^{\wedge}}(\underline{h}_{Rx},\underline{h}_{Ry})$ which in turn is weakly equivalent to
$\underline{Hom}_{L^{H}((C',W')^{\wedge})}(\underline{h}_{Rx},\underline{h}_{Ry})$. This shows that the morphism
of $S$-categories $L^{H}\underline{h}_{R}:L^{H}C'\rightarrow L^{H}((C',W')^{\wedge})$ is fully faithful.
To infer from this that the morphism  $L^{H}C'\rightarrow L^{H}\widehat{C}$ is likewise fully faithful, it is
enough to observe that we have a commutative diagram of $S$-categories
$$\xymatrix{
L^{H}((C',W')^{\wedge}) \ar[r] \ar[d] & L^{H}\widehat{C'} \ar[d] \\
L^{H}((C,W)^{\wedge}) \ar[r] & L^{H}\widehat{C} }$$
\noindent
in which the horizontal arrows are fully faithful (as $(C',W')^{\wedge}$ is a left Bousfield localization of
$\widehat{C'}$ and $(C,W)^{\wedge}$ is a left Bousfield localization of $\widehat{C}$) and the left vertical arrow
is an $S$-equivalence because $(C')^{f}$ equals $C$ (\cite[Prop. 4.1.6]{tv}).

This shows that the morphism $L^{H}C' \longrightarrow L^{H}\widehat{C}$ is fully faithful. One checks easily that as
$L^{H}C'$ is a pointed $S$-category, this also implies that the morphism $L^{H}C' \longrightarrow L^{H}\widehat{C}_{*}$
is also fully faithful.  \hfill $\Box$ \\

For the second half of the proof of Theorem \ref{t1}, let us consider the localization
morphism $l : C \longrightarrow L^{H}C$ and the induced functor on the model categories of
pointed simplicial presheaves
$$l^{*} : \widehat{L^{H}C}_{*} \longrightarrow \widehat{C}_{*}.$$
Recall that, by definition, the good Waldhausen category $M(L^{H}C)$ is the full subcategory of
$(\widehat{L^{H}C}_{*})^{f}$ consisting of representable objects\footnote{
We warn the reader that we are dealing here with two different notions of representable objects, one
in $\widehat{C}_{*}$ and the other one in $\widehat{L^{H}C}_{*}$. In the same way, we will not make
any difference between $\underline{h}_{x}$ as an object in $\widehat{C}_{*}$ or as an object
in $\widehat{L^{H}C}_{*}$ (this might be a bit confusing as $C$ and $L^{H}C$ have the same set of objects).
Note however that this abuse is justified by the fact that the standard properties of mapping spaces imply
that the ``simplicial''  $\underline{h}_{x}$ defined on page 10 coincides, up to equivalence, with the model
(or good Waldhausen, involving the choice of a cosimplicial resolution $\Gamma$) $\underline{h}_{x}$ defined on page 15.}

The Yoneda lemmas for pseudo-model categories (see \cite[Lemma 4.2.2]{tv}) and for $S$-categories
\cite[Prop. 2.4.2]{tv}) imply that an object $F \in \mathrm{Ho}(\widehat{C}_{*})$ (resp. $F'
\in \mathrm{Ho}(\widehat{L^{H}C}_{*})$) is representable if and only if there exists
an object $x \in C$ such that for any $G \in \mathrm{Ho}(\widehat{C}_{*})$ that sends weak equivalences in $C$ to equivalences of simplicial sets (resp.
for any $G'\in \mathrm{Ho}(\widehat{L^{H}C}_{*})$), one has a natural isomorphism
$$Hom_{\mathrm{Ho}(\widehat{C}_{*})}(F,G)\simeq \pi_{0}(G(x)_{*}) \qquad
(\mathrm{resp.} \quad Hom_{\mathrm{Ho}(\widehat{L^{H}C}_{*})}(F',G')\simeq \pi_{0}(G(x)_{*}) \; ).$$
Here, we have denoted by $G(x)_{*}$ the homotopy fiber of
$G(x) \longrightarrow G(*)$ at the distinguished point $* \in G(*)$ via the natural morphism
$* \rightarrow x$ (note that
in $L^{H}C$ the natural morphism $* \longrightarrow x$ is only uniquely defined
up to homotopy, which is however enough for our purposes).

\begin{lem}\label{l5}
Let $l^{*} : \widehat{L^{H}C}_{*} \longrightarrow \widehat{C}_{*}$
be the functor defined above. Then, an object $F \in \widehat{L^{H}C}_{*}$ is
representable if and only if its image $l^{*}F$ is representable in
$\widehat{C}_{*}$.
\end{lem}

\textit {Proof:} We consider the induced functor on the level of homotopy categories
$$l^{*} : \mathrm{Ho}(\widehat{L^{H}C}_{*}) \longrightarrow \mathrm{Ho}(\widehat{C}_{*}).$$
By \cite[Thm. 2.3.5]{tv} and standard properties of the left Bousfield localization (e.g. see the discussion at the end of
\cite{tv} page 19), this functor is fully faithful and its essential image consists precisely of those
functor $C^{op} \longrightarrow SSet_{*}$ sending weak equivalences in $C$ to weak equivalences of simplicial sets.

Now, let $x\in C$ and let us show that
there exists an isomorphism in $\mathrm{Ho}(\widehat{C}_{*})$,
$l^{*}(\underline{h}_{x})\simeq \underline{h}_{x}$: this will show the only if part of the lemma.
The standard properties of mapping spaces imply that
$\underline{h}_{x} \in \mathrm{Ho}(\widehat{C}_{*})$ belongs to the essential image of the functor $l^{*}$.
Therefore, as $l^{*}$ is fully faithful, to prove that $l^{*}(\underline{h}_{x})\simeq \underline{h}_{x}$, it will be
enough to show that, for any $G \in \mathrm{Ho}(\widehat{L^{H}C}_{*})$, there exists a natural isomorphism
$$Hom_{\mathrm{Ho}(\widehat{C}_{*})}(l^{*}(\underline{h}_{x}),l^{*}(G))\simeq
Hom_{\mathrm{Ho}(\widehat{C}_{*})}(\underline{h}_{x},l^{*}(G)).$$
But, again by full-faithfulness of $l^{*}$, the Yoneda lemma for $S$-categories
\cite[Prop. 2.4.2]{tv} implies that the left hand side is naturally isomorphic to
$$Hom_{\mathrm{Ho}(\widehat{L^{H}C}_{*})}(\underline{h}_{x},G)\simeq \pi_{0}(G(x)_{*}).$$
On the other hand, the Yoneda of pseudo-model categories \cite[Lemma 4.2.2]{tv} implies for the right hand side an isomorphism
$$Hom_{\mathrm{Ho}(\widehat{C}_{*})}(\underline{h}_{x},l^{*}(G))\simeq \pi_{0}(l^{*}(G)(x)_{*}).$$
As the simplicial sets $G(x)_{*}$ and $l^{*}(G)(x)_{*}$ are clearly functorially equivalent, this
shows the first part of the lemma. \\

It remains to prove that if $F \in \mathrm{Ho}(\widehat{L^{H}C}_{*})$ is such that
$l^{*}(F)$ is representable, then $F$ is itself representable. For
this, we use what we have just proved before, i.e. that
$l^{*}(\underline{h}_{x})\simeq \underline{h}_{x}$. So, if
one has $l^{*}(F)\simeq \underline{h}_{x}$, the fact that $l^{*}$
is fully faithful implies that $F\simeq \underline{h}_{x}$. \hfill $\Box$ \\

The previous lemma implies in particular that the functor $l^{*}$ restricts to an exact functor $$l^{*} : M(L^{H}C) \longrightarrow R(C).$$

\begin{lem}\label{l6}
The above exact functor $l^{*} : M(L^{H}C) \rightarrow R(C)$ is a DK-equivalence.
\end{lem}

\textit {Proof:} By \cite{dk2} (see also \cite[Thm. 2.3.5]{tv}), we know that
the induced morphism of $S$-categories
$$L^{H}l^{*} : L^{H}\widehat{L^{H}C}_{*}\longrightarrow L^{H}\widehat{C}_{*} $$
is fully faithful. As the natural morphisms
$$L^{H}R(C) \longrightarrow L^{H}\widehat{C}_{*} \qquad
L^{H}M(L^{H}C) \longrightarrow L^{H}\widehat{L^{H}C}_{*}$$
are also fully faithful, we get, in particular, that the induced morphism of $S$-categories
$$L^{H}l^{*} : L^{H}M(L^{H}C) \longrightarrow L^{H}R(C)$$
is fully faithful. Furthermore, the ``if'' part of Lemma \ref{l5} implies that this morphism
is also essentially surjective. This
proves that the exact functor of good Waldhausen categories
$$l^{*} : M(L^{H}C) \longrightarrow R(C)$$
is a DK-equivalence. \hfill $\Box$ \\

To summarize, we have defined (lemmas \ref{l4} and \ref{l6}) a diagram of DK-equivalences between good Waldhausen categories
$$\xymatrix{C \ar[r]^-{\underline{h}} & R(C) & \ar[l]_{l^{*}} M(L^{H}C).}$$
By Proposition \ref{p1}, this induces a diagram of weak equivalences on the $K$-theory
spectra
$$\xymatrix{K(C) \ar[r]^-{K(\underline{h})} & K(R(C)) & \ar[l]_-{K(l^{*})}
K(M(L^{H}C))= K(L^{H}C).}$$
This concludes the proof of Theorem \ref{t1}. \hfill $\Box$

\begin{rmk} \emph{With some work, one might be able to check that the isomorphism $K(C)\simeq K(L^{H}C)$
in the homotopy category of spectra is functorial with respect to
$DK$-equivalences of good Waldhausen categories. It is actually functorial with respect
to exact functors, but this would require the strong functoriality property of
the construction $T \mapsto K(T)$, for $S$-categories $T$, that we choosed not to discuss in this paper.}
\end{rmk}

The most important corollary of Theorem \ref{t1} is the following one, which was our original goal.
It states that the $K$-theory spectrum of a good Waldhausen category is
completely determined, up to weak equivalences, by its simplicial (or hammock) localization.

\begin{cor}\label{c1}
If $C$ and $D$ are good Waldhausen categories, and if the
$S$-categories $L^{H}C$ and $L^{H}D$ are equivalent (i.e. are isomorphic
in the homotopy category $\mathrm{Ho}(S\mathrm{-Cat})$) then the $K$-theory spectra
$K(C)$ and $K(D)$ are isomorphic in the homotopy category of spectra.
\end{cor}

\end{section}

\begin{section}{Final comments}

\textbf{$K$-theory of Segal categories.} The definition we gave of the $K$-theory spectrum of a pointed $S$-category with fibered products
(Definition \ref{d5}) makes use of Waldhausen categories and the Waldhausen construction. In a way this
is not very satisfactory as one would like to have a definition purely in terms of $S$-categories. Such a
construction surely exists but might not be so easy to describe. A major problem is that,  by mimicking Waldhausen construction,
one would like to define, for an $S$-category $T$, a new $S$-category $S_{n}T$ classifying
strings of $(n-1)$ composable morphisms in $T$, or, in other words, an object like
$T^{I(n-1)}$. However, it is well known that the naive version of $T^{I(n-1)}$
does not give the correct answer, as for example it might not be invariant under
equivalences of $S$-categories in $T$. One way to solve this problem would be to use
\textit{weak simplicial functors and weak natural transformations} as defined
in \cite{cp}. Another, completely equivalent, solution is to use the theory of Segal categories of \cite{sh,p}.

As shown in \cite{sh,p} (for an overview of results, see also \cite[Appendix]{tv})
Segal categories behave very much like categories, and many of the
standard categorical constructions are known to have reasonable analogs. There exists for example a notion
of Segal categories of functors between Segal categories, a notion of limit and colimit and more generally
of adjunctions in Segal categories, a Yoneda lemma  \dots . These constructions could
probably be used in order to define the $K$-theory spectrum of any pointed Segal category with finite limits
in a very intrinsic way and without referring to Waldhausen construction. Roughly speaking, the construction should proceed as follows.
We start from any such Segal category $A$, and consider the simplicial Segal category
$$\begin{array}{cccc}
S_{*}A : & \Delta^{op} & \longrightarrow & \textrm{Segal Cat} \\
& [n] & \mapsto & S_{n}A:=A^{I(n-1)},
\end{array}$$
where the transitions morphisms are given by various fibered products as in Waldhausen original
construction (this diagram is probably not really a simplicial Segal category, but
only a weak form of it. In other words the functor $S_{*}A$ has itself to be understood
as a morphism from $\Delta^{op}$ to the $2$-Segal category of Segal categories, see \cite{sh}).
Then we consider the simplicial diagram of maximal sub-Segal groupoids (called
\textit{interiors} in \cite[\S 2]{sh})
$$\begin{array}{cccc}
wS_{*}A : & \Delta^{op} & \longrightarrow & \textrm{Segal Groupoids} \\
& [n] & \mapsto & wS_{n}A:=(A^{I(n-1)})^{int},
\end{array}$$
and define the $K$\textit{-theory spectrum of} $A$ to be the geometric realization of this diagram
of Segal groupoids, or in other words, to be the colimit of the functor
$wS_{*}A$ computed in the $2$-Segal category of Segal groupoids.

This construction would then give a well defined morphism
$$K : \textrm{(Segal Cat)}_{*}^{ex} \longrightarrow \mathrm{Sp},$$
from the $2$-Segal category $\textrm{(Segal Cat)}_{*}^{ex}$ of pointed Segal categories with finite limits, exact functors
and equivalences between them, to the Segal category of spectra.

This theory can also be pushed further, by introducing \textit{monoidal structures}. Indeed, there exists
a notion of monoidal Segal categories, as well as symmetric monoidal Segal categories (see \cite{to}). The previously sketched
construction could then be extended to obtain $E_{\infty}$-ring spectra from pointed Segal categories with finite
limits and with an exact symmetric monoidal structure. \\

Though there are practical reasons for having a $K$-theory functor defined on the level of
Segal categories (e.g., to develop the algebraic $K$-theory of \textit{derived geometric stacks} in the sense
of \cite{tv2}), there is also a conceptual reason for it. Indeed, Segal categories are models for
$\infty$-categories for which $i$-morphisms are \textit{invertible} for all $i>1$, and therefore
the $K$-theory spectrum of a Segal category can be viewed as the $K$-theory of
an $\infty$-category. Now, the simplicial localization $L^{H}(C,S)$ of a category
$C$ with respect to a subcategory $S$ is identified in \cite[\S 8]{sh} as the
\textit{universal Segal category obtained from $C$ by formally inverting the arrows in $S$}.
From a higher categorical point of view this means that $L^{H}(C,S)$ is a model
for the $\infty$-category formally obtained from $C$ by inverting
all morphisms in $S$. In other words, $L^{H}(C,S)$ is a model for
the $\infty$-categorical version of the usual homotopy category $S^{-1}C$, and is
therefore a kind of \textit{$\infty$-homotopy category} in a very precise sense.

Thinking in these terms, Theorem \ref{t1} says that the $K$-theory of a good Waldhausen
category, while not an invariant of its usual ($0$-truncated) homotopy category, is indeed an invariant of
its $\infty$-homotopy category. \\

\medskip

\textbf{Triangulated structures.} The reader will notice that we did not consider at all triangulated structures.
This might look
surprising as in several recent works around the theme \textit{K-theory and derived categories}
the main point was to see whether one could reconstruct or not the $K$-theory from the triangulated
derived categories (see
\cite{ne,sch,ds}). From the point of view adopted in this paper, Theorem \ref{t1} tells us that,
in order to reconstruct the $K$-theory space of $C$, one only needs the $S$-category
$L^{H}C$ \textit{and nothing more}. The reason for this
is that the triangulated structure on the homotopy category $\mathrm{Ho}(C)$, when it exists, is completely determined by
the $S$-category $L^{H}C$. Indeed, both fiber and cofiber sequences can be reconstructed
from $L^{H}C$, as well as the suspension functor.

The observation that the triangulated structure can be reconstructed from the simplicial structure has lead
to a notion of \textit{stable $S$-category} (this notion was introduced by A. Hirschowitz, C. Simpson and
the first author in order to replace the old notion of triangulated category). Very similar notions
already exist in homotopy theory, as the notion of stable model category of \cite[\S 7]{ho},
of \textit{enhanced triangulated category} of \cite{bk} (see also \cite[\S 7]{sh}), and
of \textit{stable homotopy theory} of \cite{he2}.
To be a bit more precise, a stable $S$-category is an $S$-category $T$ satisfying the following three conditions.

\begin{enumerate}

\item The $S$-category $T$ is pointed.

\item The $S$-category $T$ has fibered products and fibered co-products
(i.e. $T$ and $T^{op}$ satisfy the conditions of Definition \ref{d4}).

\item The loop space functor
$$\begin{array}{cccc}
\Omega : & \pi_{0}T & \longrightarrow & \pi_{0}T \\
 & x & \mapsto & *\times_{x}^{h}*
\end{array}$$
is an equivalence of categories.

\end{enumerate}

Here, the object $*\times_{x}^{h}*$ is a fibered product of the diagram
$\xymatrix{ \textrm{*} \ar[r] & x & \ar[l] \textrm{*}}$ in $T$, in the sense of Definition \ref{d4}.

Clearly, the simplicial localization $L^{H}M$ of any
stable model category $M$ is a stable $S$-category.
Conversely, one can show that a stable $S$-category $T$ always embeds nicely in some
$L^{H}M$, for $M$ a stable model category. The homotopy category
$\pi_{0}T$ will then be equivalent to a full sub-category
of $\mathrm{Ho}(M)$ which is stable by taking homotopy fibers. In particular,
the general framework of \cite[\S 7]{ho} will imply that the category $\pi_{0}T$ possesses a natural
triangulated structure.

Our Corollary \ref{c1} implies the following result.

\begin{cor}\label{c2}
There exist two non-equivalent stable $S$-categories $T$ and $T'$, whose
associated triangulated categories $\pi_{0}T$ and $\pi_{0}T'$ are
equivalent.
\end{cor}

\textit{Proof:} Let $M_{1}:=m\mathcal{M}(\mathbb{Z}/p^{2})$ and $M_{2}:=m\mathcal{M}(\mathbb{Z}/p[\epsilon])$ be the
two stable model categories considered in \cite{sch}. The two simplicial localizations $L^{H}M_{1}$ and
$L^{H}M_{2}$ are stable $S$-categories, which by Corollary \ref{c2} and
\cite{sch} can not be equivalent. However, it is shown in \cite{sch} that the corresponding
triangulated categories $\pi_{0}L^{H}M_{1}\simeq \mathrm{Ho}(M_{1})$ and
$\pi_{0}L^{H}M_{2}\simeq \mathrm{Ho}(M_{2})$ are indeed equivalent. \hfill $\Box$ \\

We conclude in particular that a stable $S$-category $T$ contains strictly more
information than its triangulated homotopy category $\pi_{0}T$.  \\

\medskip

\textbf{$S$-Categories and ``d\'erivateurs de Grothendieck''.}
In this work we have used the construction $M \mapsto L^{H}M$, sending a
model category $M$ to its simplicial localization $L^{H}M$ as a
substitute to the construction of the homotopy category. There exists another
natural construction associated to a model category $M$,
the \textit{d\'erivateur} $\mathbb{D}(M)$ of $M$, which was introduced
by A. Heller in \cite{he} and
by A. Grothendieck in \cite{gr} (see \cite{mal} for more detailed references).
The object $\mathbb{D}(M)$ consists essentially
of the datum of the $2$-functor sending a category $I$ to the homotopy category
$\mathrm{Ho}(M^{I})$, of $I$-diagrams in $M$.

Its seems very likely that
the strictification theorem \cite[Thm. 18.5]{sh} (see also \cite[Thm. A.3.3]{tv} or \cite[Thm. 4.2.1]{msri})
together with
the results of \cite{dk2} imply that both objects $L^{H}M$ and
$\mathbb{D}(M)$ determine more or less each others\footnote{
This has to be understood in a very weak sense. To be a bit more precise
the $S$-category $L^{H}M$ seems to reconstruct completely $\mathbb{D}(M)$, but
$\mathbb{D}(M)$ only determines $L^{H}M$ as an object in the homotopy category
of $S$-categories $\mathrm{Ho}(S-Cat)$. In particular some higher homotopical
information is lost when passing from $L^{H}M$ to $\mathbb{D}(M)$.
For example the simplicial monoid of self-equivalences of $L^{H}M$ seems out
of reach from $\mathbb{D}(M)$.}
and therefore should capture
roughly the same kind of homotopical information from $M$. One should be able to check for example
that for two model categories $M$ and $M'$, $L^{H}M$ and $L^{H}M'$ are equivalent if
and only if $\mathbb{D}(M)$ and $\mathbb{D}(M')$ are equivalent. Therefore, our
reconstruction theorem \ref{t1} suggests that the $K$-theory of a reasonable Waldhausen category
is more or less an invariant of its associated d\'erivateur, and there have already been some conjectures
in this direction by G. Maltsiniotis (see \cite{mal2}).

However, we would like to mention
that the obvious generalization of conjecture $1$ of \cite{mal2} to all Waldhausen
categories can not be true for obvious functoriality reasons. Indeed,
if true for all Waldhausen categories, \cite[Conjecture $1$]{mal2} would imply that
the Waldhausen $K$-theory of spaces $X \mapsto K(X)$ would factor, up to a natural equivalence,
through the category of pr\'e-d\'erivateurs (the category of $2$-functors from
$\mathrm{Cat}^{op}$ to $\mathrm{Cat}$). This would imply that the natural morphism induced on the simplicial monoids
of self-equivalences
$$\underline{Aut}(X) \longrightarrow \underline{Aut}(K(X))$$
would factor through the simplicial monoid of self equivalences of some object
in the category of pr\'e-d\'erivateurs, which is easily seen to be $1$-truncated (it is
equivalent to the nerve of the category of self natural equivalences and isomorphisms
between them). Therefore, \cite[Conjecture $1$]{mal2} would imply that the morphism
$$\underline{Aut}(X) \longrightarrow \underline{Aut}(K(X))$$
factors through the $1$-truncation of $\underline{Aut}(X)$. But, this is clearly false as
$K(X)$ contains the stabilization $\Omega^{\infty}S^{\infty}(X)$ as a direct factor, and
the action of $\underline{Aut}(X)$ on $\Omega^{\infty}S^{\infty}(X)$ is not
$1$-truncated for a general $X$.

We think that this observation, though not strictly speaking a counter-example to the original conjecture
$1$ on \cite{mal2}, suggests that there is no reasonable way to define
a $K$-theory functor directly on the level of Grothendieck d\'erivateurs, in the same
way as there is no reasonable $K$-theory functor defined on the level
of triangulated categories. \\

\medskip

\textbf{$S$-Categories and derived equivalences.}
Recently, D. Dugger and B. Shipley have shown that
if two rings $k$ and $k'$ have equivalent triangulated derived categories
then their $K$-theory spectra $K(k)$ and $K(k')$ are equivalent (see \cite{ds}). We would like
to mention here that our reconstruction theorem \ref{t1} and its main corollary
\ref{c1} are results of different nature and can not be recovered by
the techniques of \cite{ds}. Indeed, in \cite{ds} the authors only deal with very particular
type of Waldhausen categories, the categories of complexes over some
rings, which from a homotopical point of view behave in a very
particular way (see Remarks $2.5$ and $6.8$ of \cite{ds}). For example,
our results allows one to reconstruct the $K$-theory spectra of
some ring spectra $R$, whereas the techniques of \cite{ds}
do not apply in this case (in fact, there are examples of two ring spectra
with equivalent triangulated homotopy categories of modules but with
non-equivalent $S$-categories of modules). In some sense the results of the present
paper explain Remarks $2.5$ and $6.8$ of \cite{ds}, and show that
the only missing information in order to reconstruct the $K$-theory spectrum
of some Waldhausen category from its triangulated homotopy category is
encoded the mapping spaces and their composition. From our point of view, the triangulated
structure is a way to catch a bit of this information, but in general it does not see all of it. \\

\end{section}

\bigskip
\bigskip

\textbf{Acknowledgments.} Originally, the feeling that Corollary \ref{unotre} might have been true came after a question asked a few years ago by P. Bressler, concerning the possibility of defining
$K$-theory, cyclic cohomology and the Chern character directly on the level of \textit{Segal categories}.
Since then, the main ideas of this paper have been circulating informally, and we would like to thank
J. Borger, A. Neeman, M. Schlichting and C. Simpson for conversations that have motivated us to write up a detailed proof. \\
The second-named author was partially supported by the University of Bologna, funds
for selected research topics.


\medskip

\begin{thebibliography}{90}

\bibitem{bk} A. Bondal, M. Kapranov, Enhanced triangulated categories, Math. USSR-Sbornik, Vol. 70 (1991), 93--107.

\bibitem{cis} D. Cisinski, $K$-th\'eorie de Waldhausen et $K$-th\'eorie d\'eriv\'ee,
preprint available at http://www.math.jussieu.fr/$^{\sim}$cisinski.

\bibitem{cp} J.M. Cordier, T. Porter, Homotopy coherent category theory,
Trans. Amer. Math. Soc. 349, 1-54.

\bibitem{ds} D. Dugger, B. Shipley, $K$-theory and derived equivalences,
Preprint math.KT/0209084.

\bibitem{dk1} W. Dwyer,  D. Kan, Simplicial localization of categories,
J. Pure and Appl. Algebra 17 (1980), 267-284.

\bibitem{dk2} W. Dwyer,  D. Kan, Equivalences between homotopy theories of diagrams,
in: Algebraic topology and algebraic $K$-theory,
Annals of Math. Studies 113, Princeton
University Press, Princeton, 1987, 180--205.

\bibitem{dk3} W. Dwyer,  D. Kan, Calculating simplicial localizations,
J. Pure and Appl. Algebra 18 (1980), 17--35.

\bibitem{dk4} W. Dwyer,  D. Kan, Function complexes in homotopical algebra,
Topology 19 (1980), 427--440.

\bibitem{gr} A. Grothendieck, D\'erivateurs, unpublished manuscript, around 1990.

\bibitem{he} A. Heller, \textit{Homotopy theories}, Memoirs of the Amer. Math. Soc. 71
(383), 1988.

\bibitem{he2} A. Heller, Stable homotopy theories and stabilization, J. Pure Appl. Algebra 115
(1997), 113--130.

\bibitem{hi} P. S. Hirschhorn, Model Categories and Their Localizations, Mathematical Surveys and Monographs, Vol. 99, AMS, Providence, 2003.

\bibitem{sh} A. Hirschowitz, C. Simpson, Descente pour les $n$-champs,
Preprint math.AG/9807049.

\bibitem{ho} M. Hovey, Model categories, Mathematical surveys and monographs, Vol. 63,
Amer. Math. Soc. Providence 1998.

\bibitem{ho2} M. Hovey, Model categories structures on chain complexes of sheaves, Transactions of the AMS 353 (2001), 2441--2457.

\bibitem{mal} G. Maltsiniotis, Introduction \`a la th\'eorie des d\'erivateurs, preprint available at http://www.math.jussieu.fr/$^{\sim}$maltsin.

\bibitem{mal2} G. Maltsiniotis, La $K$-th\'eorie d'un d\'erivateur triangul\'e, preprint available at http://www.math.jussieu.fr/$^{\sim}$maltsin.

\bibitem{ne} A. Neeman, $K$-Theory for triangulated categories $3$ $\frac{3}{4}$, $K$-theory, 22 (2001),
1-144.

\bibitem{p} R. Pellissier, Cat\'egories enrichies faibles, Th\`ese,  Universit\'e de Nice-Sophia Antipolis,
June 2002, Preprint math.AT/0308246.

\bibitem{sch} M. Schlichting, A note on K-theory and triangulated categories, Inv. Math. 150 (2002), 111--116.

\bibitem{tt} R. W. Thomason, T. Trobaugh, Higher algebraic $K$-theory of schemes and of derived categories,  in: The Grothendieck Festschrift, Vol. III, Birkh\"auser, 1990, pp. 247-435.

\bibitem{to} B. To\"en, Dualit\'e de Tannaka sup\'erieure $I$: Structures mono\"{\i}dales, MPI Preprint, available at http://www.mpim-bonn.mpg.de.


\bibitem{tv} B. To\"en, G. Vezzosi, Homotopy Algebraic Geometry I:
Topos theory, preprint math.AG/0207028, submitted.

\bibitem{tv2} B. To\"en, G. Vezzosi, Homotopy Algebraic Geometry II: Derived geometric stacks, in preparation.

\bibitem{msri} B. To\"en, G. Vezzosi, Segal topoi and stacks over Segal sites, to appear in: Proceedings of the program Stacks, Intersection theory and Non-abelian Hodge Theory, MSRI  Berkeley, January-May 2002.

\bibitem{wa} F. Waldhausen, Algebraic $K$-theory of spaces, in:
A. Ranicki, N. Levitt, F. Quinn (eds.), Algebraic and Geometric Topology, Lecture Notes in Mathematics \textbf{1126},
Springer $1985$, pp. 318-419.

\end{thebibliography}
\end{document}